\documentclass[11pt]{article}
\usepackage{amsfonts}
\usepackage{amssymb}
\usepackage{amsthm}
\usepackage{amsmath}
\usepackage{latexsym}
\usepackage[T1]{fontenc}
\usepackage[latin1]{inputenc}
\usepackage{graphicx}
\usepackage{amssymb}
\usepackage{amsthm}
\usepackage{amsmath}
\usepackage{latexsym}
\usepackage{color}

\newtheorem{Theo}{Theorem}
\newtheorem{Lem}{Lemma}
\newtheorem{Prop}{Proposition}
\newtheorem{Rem}{Remark}

\newtheorem{Def}{Definition}
\newtheorem{Cor}{Corollary}

\begin{document}
\title{Quadratic BSDEs with rough drivers and $L^2$--terminal condition}

\author{\textsf{M'hamed Eddahbi}$^{\mathrm{a}%
} $ and \textsf{Abou Sène}$^{\mathrm{b}}$ \\
%EndAName
$^{\mathrm{a}}${\small \textit{Cadi Ayyad University, Faculty of Sciences and Techniques, }}\\
{\small \textit{Department of Mathematics, B.P. 549, Marrakech, Morocco.}}\\
{\small \ e--mail: m.eddahbi@uca.ma} \\
$^{\mathrm{b}}${\small \textit{Gaston Berger University, UFR of Applied Sciences
and Technology,}}\\
{\small \ \textit{Department of Mathematics, B.P. 234, Saint--Louis,
Senegal}
}\\
{\small \ e--mail: seneugb@yahoo.fr}}
\maketitle

\begin{abstract}
In this paper, we study the existence and uniqueness of solutions to
quadratic Backward Stochastic Differential Equations (QBSDEs for short) with
rough driver and square integrable terminal condition. The main idea
consists in using both Doss-Sussman and Zvonkin type transformations. As an
application we study connection between QBSDEs and quadratic PDEs with rough
drivers. We also obtain Backward Doubly SDEs and QBSDEs driven by Fractional
Brownian with Hurst parameter greater than $\frac{1}{4}$ as particular cases
of our QBSDEs with rough drivers. A probabilistic representation of a class
of rough quadratic PDE is also proved.
\end{abstract}

\noindent \textbf{Keywords}: Quadratic backward stochastic differential equations, rough paths,
Doss-Sussman transformation, Zvonkin type transformation, Itô--Krylov
formula, quadratic rough PDE.

\markboth{M. Eddahbi and A. S\`ene}{Quadratic BSDEs with rough drivers and
applications}

%\ccode{AMS Subject Classification: 60H10, 60H15}

\section{Introduction}

On a Brownian motion setting a Backward stochastic differential equations
(BSDE) with terminal variable $\xi$ at time horizon $T$ and generator $g$ is
solved by a pair of processes $(Y, Z)$ on the interval $[0, T ]$ satisfying
\[
Y_{t}=\xi +\int_{t}^{T}g(s,Y_{s},Z_{s})ds-\int_{t}^{T}Z_{s}dW_{s}, \ 0\leq
t\leq T,
\]
where $(W_{t})_{0\leq t\leq T}$ is a standard $d$--dimensional Brownian
motion. Due to their interesting applications in control theory and in
partial differential equations (PDE) they have been extensively studied
since the first paper of Pardoux and Peng \cite{papeng} where they proved
that there exists a unique solution to this equation when the terminal
condition $\xi $ and the coefficient $g$ satisfy smooth square integrability
assumptions and if $g(t,\omega ,y,z)$ is Lipschitz in $(y,z)$ uniformly in $%
(t,\omega )$. Since then, several contributions have been done for relaxing
those assumptions. Kobylanski \cite{koby} studied a BSDE when the generator $%
g$ is continuous and has a quadratic growth in $z$ and the terminal
condition is bounded. Since then, there were many works on QBSDE. We notice
that, all established results on QBSDE require the terminal condition to be
bounded or of finite exponential moments, see among others \cite{brian},
\cite{lepele}. Recently, Bahlali \textit{et al}. \cite{kaled} have studied
one-dimensional QBSDE with square integrable terminal value. More precisely
they established existence and uniqueness of square integrable solutions for
a class of QBSDE when the generator $g$ is dominated by a generator of the
form $f(y)|z|^2$ where $f$ is measurable and integrable over $\mathbb{R}$.

Motivated by rough path PDE (see for instance \cite{car}, \cite{cara}, \cite{deya} and \cite{gubi}). Diehl and Friz \cite{joch} considered
\begin{equation}
Y_{t}=\xi +\int_{t}^{T}g(s,Y_{s},Z_{s})ds+\int_{t}^{T}G(Y_{s})d\mathbf{\eta}
_{s}-\int_{t}^{T}Z_{s}dW_{s},~~0\leq t\leq T,  \label{eq1}
\end{equation}
where $G=(G_{1},\ldots ,G_{d})$ is a vector field in $\mathbb{R}$, with $%
G_{k}:~\mathbb{R}~\longrightarrow ~\mathbb{R},~k=1,\dots ,d$. $\mathbf{\eta}
$ is a general geometric rough path, which by definition means that there
exists a sequence of smooth paths ($\eta^{n}$) converging to $\mathbf{\eta}$
in $p$--variation rough path metric, $p\geq 1$. When $\xi $ is bounded and
the function $g$ is of quadratic growth in $z$, they used stability theory
developed in \cite{koby} and proved existence and uniqueness of solutions to
QBSDE (\ref{eq1}). The theory of rough paths has been the subject of several
papers and lecture notes. We refer for instance to (\cite{friz}, \cite{ley},
\cite{tery}, \cite{terry} and \cite{terri}) for interesting research works
in this domain.

The first aim of this paper is the study of the QBSDE (\ref{eq1}) when the
terminal data $\xi$ is square integrable and the generator $g$ satisfy a
quadratic growth condition in $z$ to be specified the assumptions below.
Under some hypothesis weaker than those in \cite{joch} we prove existence of
a solution to the QBSDE (\ref{eq1}). In some particular cases we establish
also uniqueness. Indeed the boundedness of the random variable $\xi $, the
continuity of the generator $g$ and the linear growth of its partial
derivatives are not needed for the existence and uniqueness of solutions.
The main tools is to use Doss-Sussman and Zvonkin transformations.

The second aim consists in giving probabilistic representations of some
quadratic PDE with rough drivers in the Markovian framework. Moreover some
particular rough path are presented to deduce form our first result
existence and/or uniqueness for backward doubly stochastic differential
equations (BDSDE) and QBSDE driven by fractional Brownian motion (fBm) with
Hurst parameter greater than $\frac{1}{4}\cdot$

Our contribution presents clearly as generalization of the main results of
Bahlali \textit{et al}. \cite{kaled} and \cite{BEO} to QBSDE with rough
drivers and extend the results of \cite{joch} to merely square integrable
terminal data and measurable and integrable generator.

The rest of this paper is organized as follows. Section $2$ is devoted to
the notations, definitions and the main assumptions on the data. In Section~$%
3$ we state and prove our main result concerning the existence of QBSDE with
rough driver and square integrable terminal condition $\xi $. In Section $4$
we study the solvability of a class of QBSDE with rough driver by using a
Zvonkin transformation. In Section $5$ we restrict our selves to a Markovian
setting and study a class of PDE driven by rough path. In Section $6$ we
establish the connection to BDSDE. In this context, firstly $\mathbf{\eta }$
is replaced by a standard Brownian motion. Secondly we investigate the case
when $\mathbf{\eta }$ is replaced by a fBm with Hurst parameter $H>\frac{1}{4%
}\cdot $

\section{Notations, definitions and assumptions}

We fix once and for all a time interval $[0, T]$ and a filtered probability
space $(\Omega ,\mathcal{F},(\mathcal{F}_{t})_{0\leq t \leq T},\mathbb{P})$
which carries a $d $--dimensional Brownian motion $W$. For a vector $x$ we
denote the euclidean norm as usual by $|x|$. Let $(\mathcal{F}_{t})_{0\leq t
\leq T}$ be the usual filtration of $W$. In order to simplify the notations,
we sometimes write $Y$ for the process $(Y_{t})_{0\leq t \leq T}$. Denote by
$\mathcal{M}^{2}$ the space of predictable processes $Z$ in $\mathbb{R}^{d}$
such that
\[
\left\Vert Z\right\Vert ^{2}:=\mathbb{E}\left[ \int_{0}^{T}\left\vert
Z_{s}\right\vert ^{2}ds\right] <\infty .
\]
Denote by $\mathcal{S}^{2}$ the space of $\mathbb{R}$--valued predictable
processes $Y$ such that
\[
\mathbb{E}\left[\sup_{0\leq t\leq T}|Y_{t}|^{2}\right]<\infty .
\]
\[
\mathcal{W}^{2}_{p,loc}(\mathbb{R}) := \left\{u : \mathbb{R} \rightarrow
\mathbb{R}: ~~ u, u^{\prime}, u^{\prime \prime} \in L^{p}_{loc}(\mathbb{R}%
)\right\}.
\]
\[
\mathcal{L}^{2} :=\left\{Z, ~\mathcal{F}_{t}\text{--adapted such that}
\int_{0}^{T}|Z_{s}|^{2}ds < \infty ~~ \text{a.s.}\right\}.
\]
For a matrix $M$ we denote by $|M|$, depending on the situation, either the $%
1$--norm, the $2$--norm or the $\infty $--norm. For $p\geq 1, G^{[p]}(%
\mathbb{R}^{d})$ is the free step--$[p]$ nilpotent group over the space $%
\mathbb{R}^{d}$, realized as subset of $\mathbb{R}\oplus \mathbb{R}%
^{d}\oplus (\mathbb{R}^{d})^{\otimes 2}\oplus \ldots \oplus (\mathbb{R}%
^{d})^{\otimes[p]} $, equipped with Carnot--Caratheodory norm as defined in
\cite{joch}. Here $[p]$ denotes the largest integer not larger than $p$. $%
\mathcal{C}^{p-var}_{0}([0,T],G^{[p]}(\mathbb{R}^{d}))$ is the set of
geometric $p$--variation rough paths
\[
\mathbf{\eta} :[0,T]\longrightarrow G^{[p]}(\mathbb{R}^{d}),
\]
starting from $0$. For more and technical details on geometric rough path
spaces, we refer to \cite{friz} Section $9$, but they are not necessary for
the understanding of this paper. Let $V$ be a vector space and $(\mathcal{X}%
,|~.~|)$ a Banach vector space on $\mathbb{R}^{d}$. For $\gamma >0$, we set $%
[\gamma ]:=\gamma -1$, if $\gamma \in \mathbb{N}$ and $[\gamma
]:=\left\lfloor \gamma \right\rfloor $, the integer part of $\gamma $, if $%
\gamma \notin \mathbb{N}$. For $j \in \mathbb{N}$, $d^{j}f$ denotes the
derivative of order $j$ of the differentiable function $f$. For $x \in
\mathbb{R}^{n}$, we denote $D := D_{x} = (\frac{\partial}{\partial{x_{1}}},
\ldots, \frac{\partial}{\partial{x_{n}}})$ and $D^{2} := D_{xx} = (\frac{%
\partial^{2}}{\partial{x_{i}\partial x_{j}}})_{i, j = 1}^{n}$. Referring to
\cite{friz} we recall the following definition

\begin{Def}
For $\gamma >0,f$ is a $\gamma $--Lipschitz function on $V$ if

\begin{enumerate}
\item $f:V\longrightarrow \mathcal{X}$ is $[\gamma ]$--times differentiable.

\item {$d^{j}f$ is bounded by $K$, for all $j=0,\ldots ,[\gamma ]$.}

\item {$d^{[\gamma ]}f$ is $(\gamma -[\gamma ])$--Hölder, with Hölder
constant $K$, i.e. }
\begin{equation}
\text{for all $x\neq y\in V$,}~~\frac{\left\vert d^{[\gamma
]}f(x)-d^{[\gamma ]}f(y)\right\vert }{\left\vert x-y\right\vert ^{\gamma
-[\gamma ]}}\leq K.  \label{gama}
\end{equation}
\end{enumerate}
\end{Def}

We denote by $Lip^{\gamma }(V)$ the set of $\gamma $--Lipschitz functions on
$V$. The smallest constant $K$ for which the inequality (\ref{gama}) is
satisfied is called the Lipschitz norm of $f$ and is denoted by $\left\Vert
f\right\Vert _{Lip^{\gamma }(V)}$.

In what follows, we will refer to equation (\ref{eq1}) as: BSDE$(\xi, g, G,
\mathbf{\eta})$.

\begin{Def}
We call $(Y,Z)$ a solution of the BSDE$(\xi, g, G, \mathbf{\eta})$ if

\begin{enumerate}
\item[(i)] {$(Y,Z)\in \mathcal{S}^{2}\times \mathcal{M}^{2}$.}

\item[(ii)] {For each $t\in[0,T],(Y_{t},Z_{t})$ satisfies (\ref{eq1})}.
\end{enumerate}
\end{Def}

\textbf{Assumptions:} \newline
(H1) The function $g$ is continuous in $(y,z)$, for a.e. $(t,\omega )$.
\newline
(H2) $\mathbb{P}$--a.s. for $(t,y,z)\in \lbrack 0,T]\times \mathbb{R}\times
\mathbb{R}^{m}$,
\[
\left\vert g(t,y,z)\right\vert \leq a+b|y|+c|z|+f(|y|)|z|^{2},
\]%
where $a$, $b$ and $c$ are positive real numbers (which may change from line
to line), and $f$ a positive integrable function. \newline
(H3) For given real numbers $\gamma >p\geq 1$ and $C_{G}>0$, we have
\[
|G|_{Lip^{\gamma +2}(\mathbb{R})}:=\sup_{i=1,\dots ,d}|G_{i}|_{Lip^{\gamma
+2}(\mathbb{R})}\leq C_{G}.
\]%
As a consequence of Theorem $3.1$ in \cite{kaled}, we get the following
Proposition.

\begin{Prop}
\label{propo1} Let $G$ be Lipschitz on $\mathbb{R}$, $\xi \in L^{2}(\mathcal{%
F}_{T})$ and $\eta$ a given smooth path. Under the assumptions (H1)--(H2)
the QBSDE
\begin{equation}  \label{eq3}
Y_{t}=\xi +\int_{t}^{T}g(s,Y_{s},Z_{s})ds+\int_{t}^{T}G(Y_{s})d{\eta}%
_{s}-\int_{t}^{T}Z_{s}dW_{s},~~0\leq t\leq T,
\end{equation}
has at least one solution.
\end{Prop}

\noindent\textbf{Proof.} The path $\eta $ is smooth, then we can rewrite the
QBSDE as
\begin{eqnarray*}
Y_{t} &=&\xi +\int_{t}^{T}g(s,Y_{s},Z_{s})ds+\int_{t}^{T}G(Y_{s})\dot{\eta}
_{s}ds-\int_{t}^{T}Z_{s}dW_{s} \\
&=&\xi +\int_{t}^{T}\left( g(s,Y_{s},Z_{s})+G(Y_{s})\dot{\eta} _{s}\right)
ds-\int_{t}^{T}Z_{s}dW_{s}.
\end{eqnarray*}
Since
\[
\left\vert g(s,y,z)+G(y)\dot{\eta} _{s}\right\vert \leq \left\vert
g(s,y,z)\right\vert +\left\vert G(y)\dot{\eta} _{s}\right\vert ,
\]
by (H2) and the Lipschitz property of $G$ we get
\[
\left\vert g(s,y,z)+G(y)\dot{\eta} _{s}\right\vert \leq
a+b|y|+c|z|+f(|y|)\left\vert z\right\vert ^{2}+\left\vert G(0)\right\vert
\left\vert \dot{\eta} _{s}\right\vert +C|y|\left\vert \dot{\eta}
_{s}\right\vert ,
\]
where $C$ is the Lipschitz constant of $G$. Without loss of generality, we
may assume $\dot{\eta}$ to be bounded then we obtain
\[
\left\vert g(s,y,z)+G(y)\dot{\eta} _{s}\right\vert \leq
a+b|y|+c|z|+f(|y|)\left\vert z\right\vert ^{2}.
\]
The rest of the proof follows from Theorem $3.1$ in \cite{kaled}.

\section{Main results}

\subsection{Doss-Sussman transformation}

The first main results of this paper is the existence of solutions to the
rough quadratic equation BSDE$(\xi ,g,G,\mathbf{\eta})$ when $\mathbf{\eta}
\in \mathcal{C}^{p-var}_{0}([0,T],G^{[p]}(\mathbb{R}^{d}))$) and $\xi$ is
square integrable.

\begin{Theo}
\label{mt} Let $\gamma >p\geq 1,~ \mathbf{\eta} \in \mathcal{C}%
^{p-var}_{0}([0,T],G^{[p]}(\mathbb{R}^{d}))$ and $\xi \in L^{2}(\mathcal{F}%
_{T})$. Assume (H1)-(H3) hold. Then the quadratic BSDE$(\xi ,g,G,\mathbf{\eta%
})$ has at least one solution.
\end{Theo}

To prove this theorem we need some technical Lemmas.

Let $\phi$ be the flow solution of the following Ordinary Differential
Equations (ODE)
\begin{equation}  \label{flo1}
\phi (t,y)=y+\int_{t}^{T}\sum_{k=1}^{d}G_{k}(\phi (s,y))d\eta _{s}^{k},
\end{equation}
with $y$--inverse $\phi ^{-1}$ is given by
\[
\phi ^{-1}(t,y)=y-\int_{t}^{T}\sum_{k=1}^{d}\partial _{y}\phi
^{-1}(s,y)G_{k}(y)d\eta _{s}^{k}.
\]
The following lemma gives a way to prove and construct a solution to QBSDE (%
\ref{eq3}) using Doss-Sussmann transformation to remove the term containing
the rough path.

\begin{Lem}
\label{lemme1} Let us given a Lipschitz function $G$ on $\mathbb{R}$, a
square integrable random variable $\xi $ and a smooth path $\eta $. Let $%
\phi $ be the flow defined in (\ref{flo1}). A couple $(Y,Z)$ is a solution
of the QBSDE (\ref{eq3}) if and only if the process $(\tilde{Y},\tilde{Z})$
defined as
\[
\tilde{Y}_{t}:=\phi ^{-1}(t,Y_{t}),~~\tilde{Z}_{t}:=\frac{1}{\partial
_{y}\phi (t,\tilde{Y_{t}})}Z_{t},
\]
satisfies the BSDE
\begin{equation}
\tilde{Y}_{t}=\xi +\int_{t}^{T}\tilde{g}(s,\tilde{Y_{s}},\tilde{Z_{s}}
)ds-\int_{t}^{T}\tilde{Z_{s}}dW_{s},  \label{ft}
\end{equation}
where (throughout, $\phi $ and all its derivatives will always be evaluated
at $(t,\tilde{y})$)
\[
\tilde{g}(t,\tilde{y},\tilde{z}):=\frac{1}{\partial _{y}\phi }\big(g(t,\phi
,\partial _{y}\phi \tilde{z})+\frac{1}{2}\partial _{y}^{2}\phi |\tilde{z}%
|^{2}\big).
\]
\end{Lem}

\noindent\textbf{Proof.} Suppose that the process $(Y, Z)$ is a solution of
the QBSDE (\ref{eq3}). Denoting by $\psi :=\phi ^{-1}$ and $\theta
_{s}=(s,Y_{s})$, we have by Itô's formula
\begin{eqnarray*}
\psi (\theta _{t}) &=&\xi -\int_{t}^{T}\sum_{k=1}^{d}\partial _{y}\psi
(\theta _{s})G^{k}(Y_{s})\dot{\eta} _{s}^{k}ds+\int_{t}^{T}\partial _{y}\psi
(\theta _{s})g(s,Y_{s},Z_{s})ds \\
&& \\
&&+\int_{t}^{T}\sum_{k=1}^{d}\partial _{y}\psi (\theta _{s})G^{k}(Y_{s})\dot{%
\eta}_{s}^{k}ds-\int_{t}^{T}\partial _{y}\psi (\theta _{s})Z_{s}dW_{s} \\
&& - \frac{1}{2}\int_{t}^{T}\partial _{y}^{2}\psi (\theta _{s})|Z_{s}|^{2}ds
\\
&& \\
&=&\xi +\int_{t}^{T}\left(\partial _{y}\psi (\theta _{s})g(s,Y_{s},Z_{s})-
\frac{1}{2}\partial _{y}^{2}\psi (\theta _{s})|Z_{s}|^{2}\right) ds \\
&& -\int_{t}^{T}\left\langle \partial _{y}\psi (\theta
_{s})Z_{s},dW_{s}\right\rangle .
\end{eqnarray*}
Now, by deriving the identity $\psi (t,\phi (t,\tilde{y}))=\tilde{y}$ we get
\[
1=\partial _{y}\psi \partial _{y}\phi ,~~\text{and}~~0=\partial _{y}^{2}\psi
(\partial _{y}\phi )^{2}+\partial _{y}\psi \partial _{y}^{2}\phi .
\]
And hence,
\[
\partial _{y}\psi =\frac{1}{\partial _{y}\phi }~~\text{and}~~\partial
_{y}^{2}\psi =-\frac{\partial _{y}\psi \partial _{y}^{2}\phi }{(\partial
_{y}\phi )^{2}}=-\frac{\partial _{y}^{2}\phi }{(\partial _{y}\phi )^{3}}%
\cdot
\]
Define

\[
\tilde{Y}_{t}:=\phi ^{-1}(t,Y_{t}),~~\tilde{Z}_{t}:=\frac{1}{\partial
_{y}\phi (t,\tilde{Y_{t}})}Z_{t},
\]
and ($\psi $ and its derivatives are always evaluated at $(t,\phi (t,\tilde{y%
}))$, $\phi $ and its derivatives are always evaluated at $(t,\tilde{y})$)
\[
\tilde{g}(t,\tilde{y},\tilde{z}):=\frac{1}{\partial _{y}\phi }\big(g(t,\phi
,\partial _{y}\phi \tilde{z})+\frac{1}{2}\partial _{y}^{2}\phi |\tilde{z}%
|^{2}\big).
\]

Thanks to non--explosion condition (c.f. Condition $4.3$ in \cite{friz},
Section $4$, p. $69$) for ODEs we can show that the flow $\phi $, its
derivative $\partial _{y}\phi $ and its inverse $\phi ^{-1}$ are bounded.
This also combined with a localization argument yield to the boundedness of
the map $\frac{1}{\partial _{y}\phi }\cdot $ Then we deduce that $(\tilde{Y},%
\tilde{Z})\in \mathcal{S}^{2}\times \mathcal{M}^{2}$. We therefore obtain
\begin{equation}
\tilde{Y}_{t}:=\xi +\int_{t}^{T}\tilde{g}(s,\tilde{Y_{s}},\tilde{Z_{s}}%
)ds-\int_{t}^{T}\tilde{Z_{s}}dW_{s}.  \label{fl}
\end{equation}%
To establish the converse, we reverse the transformation and apply the Itô's
formula to the process $Y_{t}=\phi (t,\tilde{Y}_{t})$. \newline
For ease of notations, we refer to equation (\ref{fl}) with data $(\xi ,%
\tilde{g},0,0)$ as BSDE$(\xi ,\tilde{g},0,0)$.

\begin{Rem}
\label{rmk1} If one takes $g(t, y, z) = f(y)z^{2}$, where $f$ is an
integrable function, then the uniqueness of solutions $(Y, Z) $ is deduced
from the uniqueness of solutions $(\tilde{Y}, \tilde{Z})$. A technical proof
will be given in Corollary $1$ when we deal with the rough path $\mathbf{\eta%
}$.
\end{Rem}

Now, we consider $\mathbf{\eta} \in \mathcal{C}^{p-var}_{0}([0,T],G^{[p]}(%
\mathbb{R}^{d})) $ and record properties (H1) and (H2) for the induced
function $\tilde{g}$ of the previous Lemma.
%On the other hand, equation (\ref{flo}) yields a flow of diffeomorphisms for the general geometric $p$--rough path $\mathbf{\eta}$.

\begin{Lem}
\label{lemme2} Let $\gamma > p\geq 1$ and $\mathbf{\eta} \in \mathcal{C}%
^{p-var}_{0}([0,T],G^{[p]}(\mathbb{R}^{d}))$. Assume (H1)-(H3) hold. Let $%
\Phi $ be the flow of the Rough Differential Equations (RDE)
\begin{equation}  \label{flo2}
\Phi (t,y)=y+\int_{t}^{T}\sum_{k=1}^{d}G_{k}(\Phi (s,y))d\mathbf{\eta}
_{s}^{k}.
\end{equation}
Then the function
\[
\tilde{g}(t,\tilde{y},\tilde{z}):=\frac{1}{\partial _{y}\Phi }\big(g(t,\Phi
,\partial _{y}\Phi \tilde{z})+\frac{1}{2}\partial _{y}^{2}\Phi |\tilde{z}%
|^{2}\big)
\]
satisfies the following properties:

\begin{enumerate}
\item[(i)] {$\tilde{g}$ is continuous.}

\item[(ii)] {There exists positive real numbers $\tilde{a},\tilde{b}$, $%
\tilde{c}$ and a positive integrable function $\tilde{f}$ such
that}
\[
\left\vert \tilde{g}(t,\tilde{y},\tilde{z})\right\vert \leq \tilde{a}+\tilde{%
b}|\tilde{y}|+\tilde{c}|\tilde{z}|+\tilde{f}(|\tilde{y}|)|\tilde{z}|^{2}.
\]
\end{enumerate}
\end{Lem}

\noindent\textbf{Proof.} \textit{(i)}. Note that the continuity of the
function $\tilde{g}$ follows from the continuity of the flow $\Phi$ and the
continuity of the function $g$.

\textit{(ii)}. By the property (H2) we have
\begin{eqnarray*}
\left\vert \tilde{g}(t,\tilde{y},\tilde{z})\right\vert &=&\left\vert \frac{1%
}{\partial _{y}\Phi }\big(g(t,\Phi ,\partial _{y}\Phi \tilde{z})+\frac{1}{2}%
\partial _{y}^{2}\Phi |\tilde{z}|^{2}\big)\right\vert \\
&& \\
&\leq &\frac{1}{|\partial _{y}\Phi |}\left( a+b|\Phi |+c|\partial _{y}\Phi ||%
\tilde{z}|\right) \\
&& \\
&&+\frac{1}{|\partial _{y}\Phi |}\left( f(|\Phi |)||\partial_{y}\Phi |^{2}|%
\tilde{z}|^{2}\right) +\frac{1}{2}\frac{|\partial _{y}^{2}\Phi |}{|\partial
_{y}\phi |}|\tilde{z}|^{2}.
\end{eqnarray*}
Set
\[
h(\tilde{y})=\frac{1}{2}\frac{\partial _{y}^{2}\Phi }{\partial _{y}\Phi }=
\frac{1}{2}\partial _{y}\log |\partial _{y}\Phi |.
\]
We use Proposition $11.11$, in \cite{friz}, Section $11$, p. $289$ to bound
the flow $\Phi$ and its derivatives. A localization argument combined with
the $p$ non--explosion condition \cite{friz}, Section $11$, definition $11.1
$, p. $282$ for RDE give the boundedness of the map $\frac{1}{\partial
_{y}\Phi }\cdot$ Hence
\[
\frac{1}{|\partial _{y}\Phi |}\left( a+b|\Phi |+c|\partial _{y}\Phi ||\tilde{%
z}|\right) \leq \tilde{a}+\tilde{b}|\tilde{z}|,
\]
and
\[
\frac{1}{|\partial _{y}\Phi |}\left( f(|\Phi |)||\partial _{y}\Phi |^{2}|%
\tilde{z}|^{2}\right) \leq \tilde{c}f(|\Phi |)|\tilde{z}|^{2}.
\]
It follows that
\begin{eqnarray*}
\tilde{g}(t,\tilde{y},\tilde{z}) &\leq & \tilde{a}+\tilde{b}|\tilde{z}|+%
\tilde{c}f(|\Phi |)|\tilde{z}|^{2}+ h(\tilde{y})|\tilde{z}|^{2} \\
&& \\
&=&\tilde{a}+\tilde{b}|\tilde{z}|+\tilde{f}(|\tilde{y}|)|\tilde{z}|^{2},
\end{eqnarray*}
where the constants $\tilde{a}$ and $\tilde{b}$ depend on $C_{G}$ and $%
\left\Vert \mathbf{\eta}\right\Vert _{p-var}$.

We are now ready to prove Theorem \ref{mt}.

\textbf{Proof of Theorem \ref{mt}.} Using properties of the function $\tilde{%
g}$ demonstrated in Lemma \ref{lemme2} and Theorem $3.1$ in \cite{kaled} see
also \cite{BEO} there exists a solution $(\tilde{Y},\tilde{Z})\in \mathcal{S}%
^{2}\times \mathcal{M} ^{2}$ to the quadratic BSDE$(\xi,\tilde{g},0,0)$.
Define
\[
Y_{t}:=\Phi (t,\tilde{Y}_{t})~~\text{ and }~~Z_{t}:=\partial _{y}\Phi (t,%
\tilde{Y}_{t})\tilde{Z}_{t},~~t\in[0,T].
\]
Since $\tilde{Y}$ is semimartingale, and $\Phi$ a rough flow of $\mathcal{C}%
^{3}$--diffeomorphism, denoting $\tilde{\theta}_{t}=(t,\tilde{Y}_{t})$ we
obtain by Itô's formula
\begin{eqnarray*}
Y_{t} &=& \! \Phi (\tilde{\theta}_{t}) \\
&=&\! \xi \! +\! \int_{t}^{T}\!\!\partial _{y}\Phi (\tilde{\theta}%
_{s})\left( \frac{1}{\partial _{y}\Phi (\tilde{\theta}_{s})}g(s,\Phi (\tilde{%
\theta} _{s}),\partial _{y}\Phi (\tilde{\theta}_{s})\tilde{Z}_{s})+\frac{1}{2%
}\partial _{y}^{2}\Phi (\tilde{\theta}_{s})|\tilde{Z}_{s}|^{2}\right)\! ds \\
&& \\
&&\! +\! \int_{t}^{T}G(\Phi (\tilde{\theta}_{s}))d\mathbf{\eta }_{s}-\frac{1%
}{2}\int_{t}^{T}\partial _{y}^{2}\Phi (\tilde{\theta}_{s})|\tilde{Z}
_{s}|^{2}ds-\int_{t}^{T}\partial _{y}\Phi (\tilde{\theta}_{s})\tilde{Z}%
_{s}dW_{s} \\
&& \\
&=&\! \xi \! +\! \int_{t}^{T}G(\Phi (\tilde{\theta}_{s}))d\mathbf{\eta }%
_{s}+\int_{t}^{T}g(s,\Phi (\tilde{\theta}_{s}),\partial _{y}\Phi (\tilde{%
\theta}_{s})\tilde{Z}_{s})ds \\
&& \\
&&\! +\frac{1}{2}\int_{t}^{T}\partial _{y}^{2}\Phi (\tilde{\theta}_{s})|%
\tilde{Z}_{s}|^{2}ds-\int_{t}^{T}\partial _{y}\Phi (\tilde{\theta}_{s})%
\tilde{Z}_{s}dW_{s}-\frac{1}{2}\int_{t}^{T}\partial _{y}^{2}\Phi (\tilde{%
\theta}_{s})|\tilde{Z}_{s}|^{2}ds \\
&& \\
&=&\! \xi \! +\! \int_{t}^{T}g(s,Y_{s},Z_{s})ds+\int_{t}^{T}G(Y_{s})d\mathbf{%
\eta }_{s}-\int_{t}^{T}Z_{s}dW_{s}.
\end{eqnarray*}
Since the rough flow $\Phi $ and its derivatives are bounded c.f.
Proposition $11.11$ in \cite{friz}, we deduce that $(Y,Z)\in \mathcal{S}%
^{2}\times \mathcal{M}^{2}$. The proof is now finished.

In the next Corollary we establish the uniqueness of solutions for a class
of rough QBSDE whose generator $g$ is of the form $f(y)|z|^{2}$ and $f$ is
merely measurable and integrable over the hole real line.

\begin{Cor}
\label{coro} Let $\gamma > p \geq 1,~ \mathbf{\eta} \in C^{p-var}_{0}([0,T],
G^{[p]}(\mathbb{R}^{d}))$ and $f$ a real--valued integrable function. Assume
that (H3) is satisfied and $\xi \in L^{2}(\mathcal{F}_{T})$. Let $\Phi$ be
the rough flow (\ref{flo2}). Then the quadratic BSDE$(\xi,~ f(y)|z|^{2},~
G,~ \mathbf{\eta})$ has a unique solution in $\mathcal{S}^{2}\times \mathcal{%
M}^{2}$.
\end{Cor}

\noindent\textbf{Proof.} Putting $g(t,y,z)=f(y)|z|^{2}$, we have
\begin{eqnarray*}
\tilde{g}(t,\tilde{y},\tilde{z}):= &&\frac{1}{\partial _{y}\Phi }\big(%
g(t,\Phi ,\partial _{y}\Phi \tilde{z})+\frac{1}{2}\partial _{y}^{2}\Phi |%
\tilde{z}|^{2}\big) \\
&& \\
&=&\frac{1}{\partial _{y}\Phi }f(\Phi )|\partial _{y}\Phi |^{2}|\tilde{z}%
|^{2}+\frac{1}{2}\frac{\partial _{y}^{2}\Phi }{\partial _{y}\Phi }|\tilde{z}%
|^{2} \\
&& \\
&=&f(\Phi )\partial _{y}\Phi |\tilde{z}|^{2}+\frac{1}{2}\partial _{y}\log
|\partial _{y}\Phi ||\tilde{z}|^{2}.
\end{eqnarray*}%
Using the boundedness of the derivative of the flow, the function $f(\Phi
)\partial _{y}\Phi $ is integrable. Moreover we can write
\[
\tilde{g}(t,\tilde{y},\tilde{z})=\tilde{f}(\tilde{y})|\tilde{z}|^{2},
\]%
where $\tilde{f}=f(\Phi )\partial _{y}\Phi +\frac{1}{2}\partial _{y}\log
|\partial _{y}\Phi |$ is an integrable function. Now, by Theorem $3.1$
(assertion (A)) in \cite{kaled} see also \cite{BEO}, we get the existence
and uniqueness of solution $(\tilde{Y},\tilde{Z})\in \mathcal{S}^{2}\times
\mathcal{M}^{2}$ for the quadratic BSDE$(\xi ,~\tilde{f}(\tilde{y})|\tilde{z}%
|^{2},0,0)$. More precisely we deduce the existence of solution $(Y,Z)\in
\mathcal{S}^{2}\times \mathcal{M}^{2}$ from the proof of Theorem \ref{mt}
for the rough quadratic BSDE$(\xi ,~f(y)|z|^{2},~G,~\mathbf{\eta })$. Its
uniqueness follows from the uniqueness of $(\tilde{Y},\tilde{Z})$ and the
fact that the mapping
\[
(Y,Z):=L(\tilde{Y},\tilde{Z}):=\left( \Phi (\cdot,\tilde{Y}%
_{\cdot}),~\partial _{y}\Phi (\cdot ,\tilde{Y}_{\cdot })\tilde{Z}_{\cdot
}\right)
\]%
is one to one.

\subsection{Zvonkin transformation}

The second result of this paper is to give another method to obtain
existence and uniqueness results for a large class of QBSDE with rough
drivers and square integrable terminal data. In many situations we can
assume that the generator of the QBSDE is merely only measurable and
integrable in $y$. Our idea is, when the generator is the form $a + b|y| +
c|z| + f(y)z^{2}$, the quadratic part of this last can eliminated by using
Zvonkin transformation which allows us to deduce the existence of QBSDE$%
(\xi, a + b|y| + c|z| + f(y)z^{2}, G, \mathbf{\eta})$ from a BSDE of the
form BSDE$(\xi, a + b|y| + c|z|, G, \mathbf{\eta})$.

We consider the following assumption and recall the It\^o--Krylov formula in
QBSDE given in (\cite{kaled}, Theorem $2.1$, \cite{BEO}).

(H4) There exists a positive stochastic process $\zeta_{t} \in L^{1}([0,
T]\times \Omega)$ and a locally integrable function $f$ such that, for every
$(t, \omega, y, z)$
\[
\left|g(t, y, z)\right| \leq \zeta_{t} + |f(y)||z|^2.
\]

\begin{Theo}
(Itô--Krylov's formula for BSDEs). \newline
Let $\xi $ be an $\mathcal{F}_{T}$--measurable and square integrable random
variable and assume that (H4) holds. Let the process $(Y,Z)$ be a solution
of the quadratic BSDE$(\xi ,g,0,0)$ in $\mathcal{S}^{2}\times \mathcal{L}%
^{2} $. Assume moreover that $\int_{0}^{T}|g(s,Y_{s},Z_{s})|ds$ is finite $%
\mathbb{P}$--almost surely. Then, for any function $u$ belonging to $%
\mathcal{C}^{1}(\mathbb{R})\cap\mathcal{W}_{1,loc}^{2}(\mathbb{R})$, we have
\[
u(Y_{t})=u(Y_{0})+\int_{0}^{t}u^{\prime }(Y_{s})dY_{s}+\frac{1}{2}%
\int_{0}^{t}u^{\prime \prime }(Y_{s})|Z_{s}|^{2}ds.
\]
\end{Theo}

The following Lemma is key element in our approach.

\begin{Lem}
\label{lemeu} Let $f$ belongs to $L^{1}(\mathbb{R})$. The function
\[
u(x): = \int_{0}^{x}\exp\left(2\int_{0}^{y}f(t)dt\right)dy
\]
satisfies the following properties

\begin{enumerate}
\item {$u^{\prime \prime }(x) - 2f(x)u^{\prime }(x) = 0 $ and $~~ u \in
\mathcal{C}^{1}(\mathbb{R})\cap \mathcal{W}^{2}_{1,loc}(\mathbb{R}).$}

\item {$u$ is a one to one function from $\mathbb{R}$ onto $\mathbb{R}$.}

\item {The inverse function $u^{-1}$ belongs to $\mathcal{C}^{1}(\mathbb{R}%
)\cap \mathcal{W}^{2}_{1,loc}(\mathbb{R}).$}

\item {Both $u$ and $u^{-1}$ are quasi--isometries.}

\item {If in addition $f$ is continuous then both $u$ and $u^{-1}$ are of $%
\mathcal{C}^{2}$.}
\end{enumerate}
\end{Lem}

\noindent\textbf{Proof.} Using the fact that $f$ is the distributional
derivative of $x \rightarrow \int_{0}^{x}f(t)dt$, we obtain the statement 1.
The $5^{th}$ one is obvious. For the rest of the proof, we refer to \cite%
{kale}, \cite{kaled}, \cite{BEO} and \cite{zon} for more details. \newline

The second result of this paper is stated in the following theorem

\begin{Theo}
Let $\gamma > p\geq 1$ and $\mathbf{\eta} \in C^{p-var}_{0}([0,T], G^{[p]}(%
\mathbb{R}^{d}))$. Assume $\xi \in L^{2}(\mathcal{F}_{T})$ and (H3) hold.
Let $f \in L^{1}(\mathbb{R})$ and $u$ the corresponding function defined in
Lemma \ref{lemeu}. Then,

\begin{enumerate}
\item[(i)] {$(Y, Z)$ is the unique solution in $\mathcal{S}^{2} \times
\mathcal{M}^{2}$ of the quadratic BSDE$(\xi, f(y)z^{2}, G, \mathbf{\eta})$
if and only if the process $(\widetilde{Y}, \widetilde{Z})$ defined as
\[
\widetilde{Y }:= u(Y), ~~ \text{and} ~~ \widetilde{Z}:= u^{\prime }(Y)Z,
\]
is the unique solution in $\mathcal{S}^{2} \times \mathcal{M}^{2}$ to the
BSDE}
\[
\widetilde{Y}_{t} = \widetilde{\xi }+ \int_{t}^{T}\widetilde{G}(\widetilde{Y
}_{s})d\mathbf{\eta}_{s} - \int_{t}^{T}\widetilde{Z}_{s}dW_{s},
\]
where $\widetilde{G }(x) = u^{\prime -1}(x)G(u^{-1}(x))$ ~ and ~~$\widetilde{%
\xi}:= u(\xi).$

\item[(ii)] {For any $a, b, c \in \mathbb{R}$, the process $(Y, Z)$ is a
solution in $\mathcal{S}^{2} \times \mathcal{M}^{2}$ of the quadratic BSDE$%
(\xi, a + b|y| + c|z| + f(y)z^{2}, G, \mathbf{\eta})$ if and only if the
process $(\widetilde{Y}, \widetilde{Z})$ is a solution in $\mathcal{S}^{2}
\times \mathcal{M}^{2}$ to the BSDE}
\begin{equation}  \label{G}
\widetilde{Y}_{t} = \widetilde{\xi}+ \int_{t}^{T}\widetilde{g}(s, \widetilde{%
Y}_{s}, \widetilde{Z}_{s})ds + \int_{t}^{T}\widetilde{G}(\widetilde{Y}_{s})d%
\mathbf{\eta}_{s} - \int_{t}^{T}\widetilde{Z}_{s}dW_{s},
\end{equation}
where $\widetilde{g}(t, y, z) = u^{\prime -1}(y)(a + b|u^{-1}(y)|)+ c|z|.$
\end{enumerate}
\end{Theo}

\noindent \textbf{Proof.} \textit{(i)}. Let $(Y,Z)$ be the unique solution
of the rough quadratic BSDE$(\xi ,f(y)z^{2},G,\mathbf{\eta })$. Since $u\in
\mathcal{C}^{1}(\mathbb{R})\cap \mathcal{W}_{1,loc}^{2}(\mathbb{R})$, then
by Itô--Krylov's formula we have
\begin{eqnarray*}
u(Y_{t}) &=&u(\xi )-\int_{t}^{T}u^{\prime }(Y_{s})dY_{s}-\frac{1}{2}%
\int_{t}^{T}u^{\prime \prime }(Y_{s})Z_{s}^{2}ds \\
&& \\
&=&u(\xi )+\int_{t}^{T}u^{\prime
}(Y_{s})f(Y_{s})Z_{s}^{2}ds+\int_{t}^{T}u^{\prime }(Y_{s})G(Y_{s})d\mathbf{%
\eta }_{s} \\
&& \\
&&-\int_{t}^{T}u^{\prime }(Y_{s})Z_{s}dW_{s}-\frac{1}{2}\int_{t}^{T}u^{%
\prime \prime }(Y_{s})Z_{s}^{2}ds \\
&& \\
&=&u(\xi )+\int_{t}^{T}\left( u^{\prime }(Y_{s})f(Y_{s})-\frac{1}{2}%
u^{\prime \prime }(Y_{s})\right) Z_{s}^{2}ds \\
&& \\
&&+\int_{t}^{T}u^{\prime }(Y_{s})G(Y_{s})d\mathbf{\eta }_{s}-\int_{t}^{T}u^{%
\prime }(Y_{s})Z_{s}dW_{s} \\
&& \\
&=&u(\xi )-\frac{1}{2}\int_{t}^{T}\left( u^{\prime \prime
}(Y_{s})-2f(Y_{s})u^{\prime }(Y_{s})\right) Z_{s}^{2}ds \\
&& \\
&&+\int_{t}^{T}u^{\prime }(Y_{s})G(Y_{s})d\mathbf{\eta }_{s}-\int_{t}^{T}u^{%
\prime }(Y_{s})Z_{s}dW_{s}.
\end{eqnarray*}%
By the assertion $1$ of the Lemma \ref{lemeu} we have
\[
u(Y_{t})=u(\xi )+\int_{t}^{T}u^{\prime }(Y_{s})G(Y_{s})d\mathbf{\eta }%
_{s}-\int_{t}^{T}u^{\prime }(Y_{s})Z_{s}dW_{s}.
\]%
We put,
\[
\widetilde{Y}_{t}:=u(Y),~~\widetilde{Z}:=u^{\prime }(Y)Z,~~\text{and}~~%
\widetilde{\xi }:=u(\xi ),
\]%
then $(\widetilde{Y},\widetilde{Z})$ satisfies the BSDE
\[
\label{king}\widetilde{Y}_{t}=\widetilde{\xi }+\int_{t}^{T}\widetilde{G}(%
\widetilde{Y}_{s})d\mathbf{\eta }_{s}-\int_{t}^{T}\widetilde{Z}_{s}dW_{s},
\]%
where $\widetilde{G}(x)=u^{\prime -1}(x)G(u^{-1}(x))$.\newline
Since
\[
|\widetilde{\xi }|\leq |\xi |\exp \left( \left\Vert f\right\Vert _{L^{1}(%
\mathbb{R})}\right) ,
\]%
then
\[
\widetilde{\xi }\in L^{2}(\mathcal{F}_{T}).
\]%
The function $u$ is Lipschitz then
\[
|\widetilde{Y}_{t}|=|u(Y_{t})|\leq C|Y_{t}|+|u(0)|,
\]%
where $C$ is a constant.\newline
Hence
\[
\widetilde{Y}\in \mathcal{S}^{2}.
\]%
By the uniform boundedness of $u^{\prime }$, we get
\[
\widetilde{Z}\in \mathcal{M}^{2}.
\]%
Since $u^{\prime }$ is uniformly bounded, the hypothesis (H3) suffices to
obtain uniqueness of solutions by Corollary \ref{coro}. Hence, the process $(%
\widetilde{Y},\widetilde{Z})$ is the unique solution in $\mathcal{S}%
^{2}\times \mathcal{M}^{2}$ of the BSDE$(\widetilde{\xi },0,\widetilde{G},%
\mathbf{\eta })$.

Conversely, suppose that the process $(\widetilde{Y},\widetilde{Z})$ is the
unique solution in $\mathcal{S}^{2}\times \mathcal{M}^{2}$ of the BSDE$(%
\widetilde{\xi },0,\widetilde{G},\mathbf{\eta })$. Since $u^{-1}\in \mathcal{%
C}^{1}(\mathbb{R})\cap \mathcal{W}_{1,loc}^{2}(\mathbb{R})$, we have again
by Itô--Krylov's formula
\begin{eqnarray*}
Y_{t} &=&u^{-1}(\widetilde{Y}_{t}) \\
&& \\
&=&u^{-1}(\widetilde{Y}_{T})-\int_{t}^{T}(u^{-1})^{\prime }(\widetilde{Y}%
_{s})d\widetilde{Y}_{s}-\frac{1}{2}\int_{t}^{T}(u^{-1})^{\prime \prime }(%
\widetilde{Y}_{s})\widetilde{Z}_{s}^{2}ds \\
&& \\
&=&\xi +\int_{t}^{T}(u^{-1})^{\prime }(\widetilde{Y}_{s})\widetilde{G}(%
\widetilde{Y}_{s})d\mathbf{\eta }_{s}-\int_{t}^{T}(u^{-1})^{\prime }(%
\widetilde{Y}_{s})\widetilde{Z}_{s}dW_{s} \\
&& \\
&&-\frac{1}{2}\int_{t}^{T}(u^{-1})^{\prime \prime }(\widetilde{Y}_{s})%
\widetilde{Z}_{s}^{2}ds.
\end{eqnarray*}%
Calculus implies
\begin{eqnarray*}
\int_{t}^{T}(u^{-1})^{\prime }(\widetilde{Y}_{s})\widetilde{G}(\widetilde{Y}%
_{s})d\mathbf{\eta }_{s} &=&\int_{t}^{T}(u^{-1})^{\prime }(\widetilde{Y}%
_{s})u^{\prime -1}(\widetilde{Y}_{s}))G(u^{-1}(\widetilde{Y}_{s}))d\eta _{s}
\\
&& \\
&=&\int_{t}^{T}(u^{-1})^{\prime }(\widetilde{Y}_{s})u^{\prime
}(Y_{s})G(Y_{s})d\eta _{s} \\
&=&\int_{t}^{T}\frac{1}{u^{\prime }(Y_{s})}u^{\prime }(Y_{s})G(Y_{s})d\eta
_{s} \\
&=&\int_{t}^{T}G(Y_{s})d\eta _{s},
\end{eqnarray*}%
\begin{eqnarray*}
\int_{t}^{T}(u^{-1})^{\prime }(\widetilde{Y}_{s})\widetilde{Z}_{s}dW_{s}
&=&\int_{t}^{T}\frac{1}{u^{\prime }(Y_{s})}u^{\prime }(Y_{s})Z_{s}dW_{s} \\
&=&\int_{t}^{T}Z_{s}dW_{s},
\end{eqnarray*}%
and also
\[
\frac{1}{2}\int_{t}^{T}(u^{-1})^{\prime \prime }(\widetilde{Y}_{s})%
\widetilde{Z}_{s}^{2}ds=\frac{1}{2}\int_{t}^{T}(u^{-1})^{\prime \prime }(%
\widetilde{Y}_{s})(u^{\prime }(Y_{s}))^{2}Z_{s}^{2}ds.
\]%
Since
\[
(u^{-1})^{\prime \prime }(\widetilde{Y}_{s})=\left( \frac{1}{u^{\prime
}(Y_{s})}\right) ^{\prime }=-2\frac{f(Y_{s})u^{\prime }(Y_{s})}{(u^{\prime
}(Y_{s}))^{2}},
\]%
then
\[
\frac{1}{2}\int_{t}^{T}(u^{-1})^{\prime \prime }(\widetilde{Y}_{s})%
\widetilde{Z}_{s}^{2}ds=\int_{t}^{T}f(Y_{s})Z_{s}^{2}ds.
\]%
Putting things together, we obtain
\[
Y_{t}=\xi +\int_{t}^{T}f(Y_{s})Z_{s}^{2}ds+\int_{t}^{T}G(Y_{s})d\mathbf{\eta
}_{s}-\int_{t}^{T}Z_{s}dW_{s}.
\]%
Moreover
\[
|Y_{t}|=|u^{-1}(\widetilde{Y}_{t})|\leq C|\widetilde{Y}_{t}|+|u^{-1}(0)|,
\]%
since the function $u^{-1}$ is Lipschitz, therefore
\[
Y\in \mathcal{S}^{2}.
\]%
Member that $Z_{t}=\frac{\widetilde{Z}_{{t}}}{u^{\prime }(u^{-1}(\widetilde{Y%
}_{t}))}$, with the inequality
\begin{eqnarray*}
\left\vert \frac{1}{u^{\prime }(x)}\right\vert &\leq &\exp \left(
-2\int_{0}^{x}f(t)dt\right) \\
&& \\
&\leq &\exp \left( 2\int_{0}^{x}|f(t)|dt\right) \\
&& \\
&\leq &\exp \left( 2\int_{0}^{|x|}|f(t)|dt\right) \\
&& \\
&\leq &\exp \left( 2\left\Vert f\right\Vert _{L^{1}(\mathbb{R})}\right) ,
\end{eqnarray*}%
one shows that $Z$ belons to $\mathcal{M}^{2}$, which means that $(Y,Z)$ is
a solution in $\mathcal{S}^{2}\times \mathcal{M}^{2}$ of the rough quadratic
BSDE$(\xi ,f(y)z^{2},G,\mathbf{\eta })$. Its uniqueness follows from the
uniqueness of $(\widetilde{Y},\widetilde{Z})$ and the fact that the mapping
\[
(Y,Z):=\mathcal{L}(\widetilde{Y},\widetilde{Z}):=\left( u^{-1}(\widetilde{Y}%
_{\cdot }),\frac{\widetilde{Z}_{\cdot }}{u^{\prime }(u^{-1}(\widetilde{Y}%
_{\cdot }))}\right)
\]%
is one to one.

\textit{(ii)}. The proof of this assertion is similar to that of \textit{(i)}%
, so the detail of the It\^o--Krylov's formula are omitted. \newline
We only need to establish the existence of solutions to BSDE (\ref{G}).
It\^o--Krylov's formula applied to the function $u$ shows that
\begin{equation}  \label{king2}
\widetilde{Y}_{t} = \widetilde{\xi}+ \int_{t}^{T}\widetilde{g}(s, \widetilde{%
Y}_{s}, \widetilde{Z}_{s})ds + \int_{t}^{T}\widetilde{G}(\widetilde{Y}%
_{s})d\eta_{s} - \int_{t}^{T}\widetilde{Z}_{s}dW_{s},
\end{equation}
where
\[
\widetilde{g}(t, y, z) = u^{\prime -1}(y)(a + b|u^{-1}(y)|)+ c|z|.
\]
The function $\widetilde{g}$ is continuous, and when we use the boundedness
of $u^{\prime }$ and the Lipschitz property of $u^{-1}$, we get
\[
|\widetilde{g}(t, y, z)| \leq a + b|y| + c|z| \leq a + b|y| + c|z| +
|f(y)||z|^{2}.
\]
By Theorem \ref{mt} the BSDE (\ref{king2}) has a solution in $\mathcal{S}%
^{2}\times \mathcal{M}^{2}$. We use the same technique developped in \textit{%
(i)} and It\^o--Krylov's formula to get existence of solutions of
%the BSDE (\ref%{G}) follows from It\^o--Krylov's formula.
BSDE$(\xi, a + b|y| + c|z| + f(y)z^{2}, G, \mathbf{\eta})$

\begin{Rem}
In contrast to \cite{joch}, our approaches cover the BSDE with linear growth
(put $f = 0$).
\end{Rem}

\subsection{Examples of application}

\subsubsection{Connection to Backward doubly SDEs}

We do the connection with the so-called backward doubly stochastic
differential equations (BDSDE) introduced by Pardoux and Peng in \cite%
{papeng94}. We recall that on $\mathcal{C}([0,T],\mathbb{R}^{d})$ there
exists a unique Borel probability measure, is known as the $d$--dimensional
Wiener measure, so that the coordinate function $B_{t}(\omega )=\omega _{t}$
defines a Brownian motion. To begin with, let $\Omega ^{1}=\mathcal{C}([0,T],%
\mathbb{R}^{d})$ and $\Omega ^{2}=\mathcal{C}([0,T],\mathbb{R}^{m})$
equipped respectively with Wiener measures $\mathbb{P}^{1}$ and $\mathbb{P}%
^{2}$. Consider $\Omega =\Omega ^{1}\times \Omega ^{2}$ on which we define
the product measure $\mathbb{P}:=\mathbb{P}^{1}\otimes \mathbb{P}^{2}$. For $%
(\omega ^{1},\omega ^{2})\in \Omega $, we define $B(\omega ^{1},\omega
^{2}):=\omega ^{1}$. Analogously, we define $W(\omega ^{1},\omega
^{2}):=\omega ^{2}$. Hence $B$ is a $d$--dimensional Brownian motion and $W$
is an independent $m$--dimensional Brownian motion. Let $\mathcal{F}_{t}:=%
\mathcal{F}_{t,T}^{B}\vee \mathcal{F}_{0,t}^{W}$, where $\mathcal{F}%
_{t,T}^{B}:=\sigma (B_{r}:r\in \lbrack t,T])\vee \mathcal{N}^{1}$, $\mathcal{%
F}_{0,t}^{W}:=\sigma (W_{r}:r\in \lbrack 0,t])\vee \mathcal{N}^{2}$ and $%
\mathcal{N}^{i}$ is set of $\mathbb{P}^{i}$--null sets, $i=1,2$ . Note that
the collection $(\mathcal{F}_{t},t\in \lbrack 0,T])$ is neither increasing
nor decreasing, and it does not constitute a filtration.

Given $\xi$ in $L^{2}(\mathcal{F}_{T})$ Pardoux and Peng \cite{papeng94}
considered the following BDSDE
\begin{equation}  \label{pp94}
Y_{t} = \xi + \int_{t}^{T}g(s, Y_{s}, Z_{s})ds + \int_{t}^{T}G(Y_{s})\circ
dB_{s} + \int_{t}^{T}Z_{s}dW_{s}, ~~ 0 \leq t \leq T.
\end{equation}
An $\mathcal{F}$--adapted process $(Y, Z)$ is called a solution of the above
BDSDE if $\mathbb{E}[\sup_{t\leq T}|Y_{t}|^{2}] < \infty$, $\mathbb{E}%
[\int_{0}^{T}|Z_{s}|^{2}ds] < \infty$ and $\mathbb{P}$--a.s. (\ref{pp94}) is
satisfied for $0 \leq t \leq T$. Under appropriate (essentially Lipschitz)
conditions on $g$ and $G$ they establish existence and uniqueness of a
solution.

Note that in \cite{papeng94} Pardoux and Peng considered the equations (\ref%
{pp94}) where the Stratonovich integral is actually a Backward It\^o
integral. But if $G$ is smooth enough, the formulations are equivalent.

We are interested now by the connection of the BDSDEs and rough drivers.
This motivates us to take $2<p<3$ and define the lift Brownian motion to a
process with values in $\mathbb{R}^{d}\oplus so(d)$, where $so(d)$ denotes the space of
anti-symmetric $d\times d$--matrices.

%Assume $B = (B^{1},
%B^{2}, \ldots, B^{d})$ is a $d$--dimensional Brownian motion. Define enhanced Brownian motion (EBM) by

\begin{Def}
(L\'evy's area). Given a $d$--dimensional Brownian motion $B = (B^{1},
B^{2}, \ldots, B^{d})$, we define the $d$--dimensional L\'evy area $A =
(A^{i,j}: i, j \in \{1, \ldots, d\})$ as the continuous process
\[
t \longrightarrow A^{i,j}_{t} = \frac{1}{2}\left(%
\int_{0}^{t}B^{i}_{s}dB^{j}_{s} - B_{s}^{j}dB_{s}^{i}\right).
\]
\end{Def}

We note that $A_{t}$ takes values in $so(d)$. In the sequel, $\exp$ denotes
the exponential map from $\mathbb{R}^{d}\oplus so(d)$ to $G^{2}(\mathbb{R}^{d})$.
%\textbf{ce n'est pas un + mais plutot $\oplus$}
Set \[G^{2}(\mathbb{R}^{d}) :=\exp(\mathbb{R}^d \oplus so(d))
= \{(1,v,\frac{1}{2}v\otimes v+A); \,  v\in \mathbb{R}^d \text {\  and \  } A \in so(d)\}.
\]
\begin{Def}
Let $B$ and $A$ denote a $d$--dimensional Brownian motion and its Lévy area
process.  The continuous $G^{2}(\mathbb{R}^{d})$--valued process $\mathbf{B}$%
, defined by
\[
\mathbf{B}_{t}:=\exp \left( B_{t}+A_{t}\right) ,~~t\geq 0
\]%
is called enhanced Brownian motion (EBM).  $\mathbf{B}$ is precisely $d$--dimensional Brownian motion enhanced with its iterated integrals in Stratonovich sense. It is in one to one correspondence with Brownian motion enhanced with Lévy's area.
%$\exp$ denotes
%the exponential map from the Lie algebra $\mathbb{R}^d \oplus so(d)$  to the group, realized inside the truncated
%tensor algebra $G^{2}(\mathbb{R}^{d})$.
\end{Def}

The EBM $\mathbf{B}$ has finite $p$%
--variation for $p\geq 2$. by setting $\mathbf{B}=0$ on a $\mathbb{P}^{1}$%
--null sets, we can say that $\mathbf{B}$ belongs to $\mathcal{C}%
_{0}^{p-var}([0,T],G^{2}(\mathbb{R}^{d}))$. The EBM can be identified as a
special case of left-invariant Brownian motion on the Lie group $G^{2}(\mathbb{%
R}^{d})$. We refer to Section $13$ in \cite{friz} for more
details.

\begin{Theo}
Let $2<p<3$ and $\gamma >p$. We assume that (H1)--(H3) hold and $\xi \in
L^{2}(\mathcal{F}_{T})$. For every $\omega ^{1}\in \Omega ^{1}$ Then the
BSDE with rough driver for all $0\leq t\leq T$
\begin{eqnarray}\label{ee}
Y_{t}^{rp}(\omega ^{1},\cdot ) &=&\xi (\cdot
)+\int_{t}^{T}g(s,Y_{s}^{rp}(\omega ^{1},\cdot ),Z_{s}^{rp}(\omega
^{1},\cdot ))ds \\
&& \nonumber \\
&&+\int_{t}^{T}G(Y_{s}^{rp}(\omega ^{1},\cdot ))d\mathbf{B}_{s}(\omega
^{1})+\int_{t}^{T}Z_{s}^{rp}(\omega ^{1},\cdot )dW_{s}(\cdot ). \nonumber
\end{eqnarray}%
has a solution. In particular, if $g(t,y,z)=f(y)|z|^{2}$, where $f$ is an
integrable function, the uniqueness holds true for this equation. Moreover by (\ref%
{pp94}) we have for $\mathbb{P}^{1}$--a.e. $\omega ^{1}$ that $\mathbb{P}^{2}
$--a.s.
\[
Y_{t}(\omega ^{1},\cdot )=Y_{t}^{rp}(\omega ^{1},\cdot ),~~t\leq T,
\]%
and
\[
Z_{t}(\omega ^{1},\cdot )=Z_{t}^{rp}(\omega ^{1},\cdot ),~~\text{$dt\otimes
\mathbb{P}^{2}$--a.s.}
\]
%\textbf{Je n'ai pas très bien compris les deux égalités ci-dessus: le couple $(Y, Z)$ est solution unique de (\ref{pp94}) pour $g(t,y,z)=f(y)|z|^{2}$  $\mathbb{P}^{1}$--a.e. $\omega ^{1}$ that $\mathbb{P}^{2}
%$--a.s. on a l'unicit\'e de la solution de (\ref{ee}. Donc pour presque tout $\omega ^{1}$ (\ref{pp94}) et (\ref{ee}) repr\'esente la m\^eme rough BSDE puisqu'on suppos\'e \textbf{B} nul sur les $\mathbb{P}^{1}$-nul sets.)}
\end{Theo}

\noindent \textbf{Proof.} As in the proof of Theorem \ref{mt}, in the BDSDEs
setting, we eliminate the integral corresponding to the Brownian motion $B$
using the stochastic flow $\phi $, defined as the unique solution of the
stochastic differential equation in the Stratonovich sense
\[
\phi (t,\omega ^{1};y)=y+\int_{t}^{T}G(\phi (s,\omega ^{1};y))\circ
dB_{s}(\omega ^{1}).
\]%
Then $\omega ^{1}$--wise, we construct the rough flow given by
\[
\Phi (t,\omega ^{1};y)=y+\int_{t}^{T}G(\Phi (s,\omega ^{1};y))d\mathbf{B}%
_{s}(\omega ^{1}).
\]%
Therefore by Theorem \ref{mt} we obtain the result.

If $g(t,y,z)=f(y)|z|^{2}$, by Corollary \ref{coro} we get uniqueness of the
solution $(Y^{rp},Z^{rp})$. By a classical result of rough path theory, we
have for every $\omega ^{1}\in \Omega ^{1}$
\[
\Phi (.,\omega ^{1};\cdot )=\phi (.,\omega ^{1};\cdot ).
\]
Hence processes $(Y,Z)$ and $(Y^{rp},Z^{rp})$ satisfy the same BSDE.
Therefore, we get the desired result by uniqueness. %
%\textbf{(*)en effet, c'est la th\'eorie classique  des rough differential equation (RDE). Pour $\omega ^{1}$ fix\'e le mouvement brownien $B$ devient d\'eterministe, parcons\'equent $\phi$ et $\Phi$ repr\'sentent la m\^eme RDE. Comme $G$ est assez smooth d'apr\`es l'hypoth\`ese $(H3)$ on l'unicit\'e de la solution}
%%%%%%%%%%%%%%%%%%%%%%%%%%%%%%%%%%%%%%%%%%%%%%%%%%%%%%%%%%%%%%%%%%%%%%%%%%%%%%%%%% %%%%%%%%%%%%%%%%%%%%%%%%%%%%%%%%%%%%%%%%%%

\subsubsection{Connection to fractional Brownian motion}

We have seen that $d$--dimensional Brownian motion $B$ can be enhanced to a
stochastic process $\mathbf{B}$ for which every realization is a geometric $p
$--rough path, $p\in (2,3)$. Recall that $B$ is a continuous, centered
Gaussian process with independent components $(B^{1},\ldots ,B^{d})$, whose
law is fully determined by its covariance function
\[
R(s,t)=\mathbb{E}(B_{s}\otimes B_{t})=diag(s\wedge t,\ldots ,s\wedge t).
\]%
We note that this covariance function $R:=R(s,t)$ has finite $1$--variation
in $2D$. More generally, consider a $d$--dimensional
continuous, centered Gaussian process with independent components $%
X=(X_{t}^{1},\ldots ,X_{t}^{d}:t\in \lbrack 0,T])$. Again its law is fully
determined by its covariance function
\[
R(s,t)=diag(\mathbb{E}(X_{s}^{1}X_{t}^{1}),\ldots ,\mathbb{E}%
(X_{s}^{d}X_{t}^{d})),~~s,t\in \lbrack 0,T].
\]%
Let $p\in (3,4)$. When the covariance function of $X$ has finite $\rho $%
--variation for some $\rho \in \lbrack 1,2)$, Friz and Victoir (\cite{friz},
Section $15$, Theorem $5.33$) have shown the existence of unique process $%
\mathbf{X}$, which lies in $\mathcal{C}_{0}^{p-var}([0,1],G^{3}(\mathbb{R}%
^{d}))$, lifting the Gaussian process $X$ for any $p>2\rho $. This $G^{3}(%
\mathbb{R}^{d})$--valued process $\mathbf{X}$ is called the enhanced
Gaussian process and sample path realizations of $\mathbf{X}$ are called
Gaussian rough paths. Theorem $5.33$ in \cite{friz} asserts in particular
that $d$--dimensional Brownian motion can be lifted to an enhanced Gaussian
process. Other example can be obtained by considering $d$--independents
copies of fractional Brownian motion (fBm) with Hurst parameter $H\in (0,1)$%
. The resulting $\mathbb{R}^{d}$--valued fBm can be lifted to an enhanced
Gaussian process provided $H>\frac{1}{4}\cdot $ Recall that a $d$%
--dimensional fBm with Hurst parameter $H\in (0,1)$ is a Gaussian process $%
B^{H}$
\[
B_{t}^{H}:=(B_{t}^{H,1},\ldots ,B_{t}^{H,d}),~~t\geq 0,
\]%
where $B^{H,1},\ldots ,B^{H,d}$ are $d$ independents centered Gaussian
processes with covariance function
\[
R(s,t)=\frac{1}{2}\left( s^{2H}+t^{2H}-|t-s|^{2H}\right) ,~~(s,t)\in \lbrack
0,+\infty \lbrack ^{2}.
\]%
Let us consider the following stochastic flow defined as the solution of the
Stochastic Differential Equation (SDEs) in the Stratonovich sense, driven by
a fBm with Hurst parameter $H>\frac{1}{4}$
\[
\phi (t,\omega ^{1};y)=y+\int_{t}^{T}G(\phi (s,\omega ^{1};y))\circ
dB_{s}^{H}(\omega ^{1}).
\]%
If $H=\frac{1}{2}$, this equation corresponds to SDEs driven by Brownian
motion in Stratonovich sense. When $H$ is greater than $\frac{1}{2}$
existence and uniqueness of the solution are obtained by Zähle \cite{zale}
and Nualart and Rascânu \cite{nualar} and the refrences therein. In the case
$H<\frac{1}{2}$, since fBm has $\alpha $--Hölder continuous sample paths for
any $\alpha <H$, it falls into the rough paths theory. When $H\neq \frac{1}{2%
}$ the fBm is neither a semimartingale nor a Markov process. Hence, a
natural application of the rough path analysis is the stochastic calculus
with respect to the fBm.

We consider $\mathcal{F}_{t}:=\mathcal{F}_{t,T}^{B^{H}}\vee \mathcal{F}%
_{0,t}^{W}$, where $\mathcal{F}_{t,T}^{B^{H}}=\sigma (B_{s}^{H}:s\in \lbrack
t,T])\vee \mathcal{N}$ and $\mathcal{N}$ is the set of $\mathbb{P}^{1}$%
--negligible sets. The canonical processes on $\Omega =\Omega ^{1}\times
\Omega ^{2}$ are defined by $B^{H}(\omega ^{1},\omega ^{2}):=\omega ^{1}$
and $W(\omega ^{1},\omega ^{2}):=\omega ^{2}$. Hence $B^{H}$ is a $d$%
--dimensional fBm.

Let $3<p<4,~H>\frac{1}{4}$ be such that $Hp>1$. In this setting Coutin and
Qian \cite{cq}, by using dyadic approximations showed the existence of a
canonical geometric $p$--rough path $\mathbf{B}^{H}$ associated to the fBm $%
B^{H}$ with Hurst parameter $H\in (\frac{1}{4},\frac{1}{2})$. Setting $%
\mathbf{B}^{H}=0$ on $\mathcal{N}$ we assume that $\mathbf{B}^{H}$ lies in $%
\mathcal{C}_{0}^{p-var}([0,T],G^{3}(\mathbb{R}^{d}))$. Consider the $\omega
^{1}$--wise fractional rough flow
\begin{equation}
\Phi (t,\omega ^{1};y)=y+\int_{t}^{T}G(\Phi (s,\omega ^{1};y))d\mathbf{B}%
_{s}^{H}(\omega ^{1}).  \label{flo3}
\end{equation}

\begin{Theo}
Let $3<p<4,~\gamma >p$ and $H>\frac{1}{4}$ be such that $Hp>1$. We assume
that (H1)--(H3) hold and $\xi \in L^{2}(\mathcal{F}_{T})$. For every $\omega
^{1}\in \Omega ^{1}$ the BSDE with rough driver for $0\leq t\leq T$
\begin{eqnarray*}
Y_{t}^{rp}(\omega ^{1},\cdot ) &=&\xi (\cdot
)+\int_{t}^{T}g(s,Y_{s}^{rp}(\omega ^{1},\cdot ),Z_{s}^{rp}(\omega
^{1},\cdot ))ds \\
&&+\int_{t}^{T}G(Y_{s}^{rp}(\omega ^{1},\cdot ))d\mathbf{B}_{s}^{H}(\omega
^{1})+\int_{t}^{T}Z_{s}^{rp}(\omega ^{1},\cdot )dW_{s}(\cdot ),
\end{eqnarray*}%
has a solution. Moreover the uniqueness holds true for $g(t,y,z)=f(y)|z|^{2}$
and $f$ is an integrable function.
\end{Theo}

\noindent \textbf{Proof.} We use the rough flow (\ref{flo3}), Theorem \ref%
{mt} and Corollary \ref{coro}.

\section{Probabilistic representation of rough PDEs}

The aim of this section is to give a probabilistic representation of the
following Quadratic PDE with rough path:
\begin{equation}
\left\{
\begin{array}{l}
du(t,x)+\left[ \mathcal{L}u(t,x)+f(u(t,x))|u_{x}(t,x)\sigma (t,x)|^{2}\right]
dt+G(u(t,x))d{\mathbf{\eta }}_{t}=0, \\
u(T,x)=\psi (x),~x\in \mathbb{R}.%
\end{array}%
\right.  \label{Rqpde}
\end{equation}%
where
\[
\mathcal{L}u(t,x):=\frac{1}{2}\left( \sigma ^{2}u_{xx}\right) (t,x)+\left(
bu_{x}\right) (t,x).
\]%
Assume first that the equation (\ref{Rqpde}) has a classical smooth
solution. Let $X$ be the unique solution to the following forward equation
\begin{equation}
X_{t}^{s,x}=x+\int_{s}^{t}b(r,X_{r}^{s,x})dr+\int_{s}^{t}\sigma
(r,X_{r}^{s,x})dW_{r},~t\in \lbrack s,T],  \label{fsde}
\end{equation}%
where the functions $\sigma $ and $b$ are given coefficients defined as
follows: $\sigma :[0,T]\times \mathbb{R}\rightarrow \mathbb{R},$ ~$%
b:[0,T]\times \mathbb{R}\rightarrow \mathbb{R},$ $f:\mathbb{R}\rightarrow
\mathbb{R},$ and $\psi :\mathbb{R}\rightarrow \mathbb{R}.$\newline
such that:\newline
(H4) $\sigma ,~b$ are uniformly Lipschitz.\newline
(H5) $\sigma ,~b$ are of linear growth.

These condition insure existence and uniqueness of the equation (\ref{fsde}).

Consider the following rough QBSDE
\begin{equation}
Y_{t}^{s,x}=\psi
(X_{T}^{s,x})+\int_{t}^{T}f(Y_{r}^{s,x})|Z_{r}^{s,x}|^{2}dr+%
\int_{t}^{T}G(Y_{r}^{s,x})d\mathbf{\eta }_{r}-\int_{t}^{T}Z_{r}^{s,x}dB_{r}
\label{Rqbsde}
\end{equation}%
where $f$, $G$ and $\psi $ are given measurable functions such that.

\noindent (H6) $f$ is continuous and integrable, $G$ is Lipschitz and $\psi $
is continuous and $|\psi (x)|\leq K(1+|x|^{p}),~~\forall ~p\geq 1$.

Applying Itô's formula to $u(T,X_{T}^{s,x})$ yields
\begin{eqnarray*}
u(T,X_{T}^{s,x}) &=&u(t,X_{t}^{s,x})+\int_{t}^{T}\left( u_{t}+\mathcal{L}%
u\right) (r,X_{r}^{s,x})dr \\
&&+\int_{t}^{T}\left( \sigma u_{x}\right) (r,X_{r}^{s,x})dW_{r} \\
&=&u(t,X_{t}^{s,x})-\int_{t}^{T}f(u(r,X_{r}^{s,x}))|\left( u_{x}\sigma
\right) (r,X_{r}^{s,x})|^{2}dr \\
&&-\int_{t}^{T}G(u(r,X_{r}^{s,x}))d\mathbf{\eta }_{r}+\int_{t}^{T}\left(
\sigma u_{x}\right) (r,X_{r}^{s,x})dW_{r}.
\end{eqnarray*}%
which means that $(u(r,X_{r}^{s,x}),\left( \sigma u_{x}\right)
(r,X_{r}^{s,x}))_{s\leq r\leq T}$ is a solution to the (\ref{Rqbsde}).

The purpose of the this section is to study the converse.

Given a solution to the (\ref{Rqbsde}) which is unique by our result in the
previous section. We shall construct a viscosity solution to the rough
quadratic PDE (\ref{Rqpde}).

When the rough path $\mathbf{\eta }$ is replaced by a smooth path $\eta $,
then (\ref{Rqpde}) takes the form
\begin{equation}
\left\{
\begin{array}{l}
\partial _{t}u(t,x)+\mathcal{L}u(t,x)+f(u(t,x))|(u_{x}\sigma
)(t,x)|^{2}+G(u(t,x))\dot{\eta}_{t}=0, \\
u(T,x)=\psi (x),~~x\in \mathbb{R},%
\end{array}%
\right.  \label{pde}
\end{equation}

We again consider the flows associated with a smooth path $\eta $ (resp.,
rough $\mathbf{\eta }$) defined in (\ref{flo1}) (resp. in (\ref{flo2})).

%\begin{equation}  \label{flow1}
%\phi(t,y) = y +\int_{t}^{T}\sum_{k=1}^{d}G_{k}(\Phi(s, y))d\eta^{k}_{s},
%\end{equation}
%
%\begin{equation}  \label{flow2}\left(\textrm{resp.,} ~~
%\Phi(t,y) = y +\int_{t}^{T}\sum_{k=1}^{d}G_{k}(\Phi(s, y))d\mathbf{\eta}^{k}_{s}\right).
%\end{equation}

\begin{Prop}
\label{propo2} Let $G$ be a Lipschitz function on $\mathbb{R}$ and $\eta$ a
given smooth path. Under the assumptions (H4)--(H6), $u(t, x):= Y^{t,x}_{t}$
is a viscosity solution to the PDE (\ref{pde}), where for every $(s, x) \in
[0, T]\times \mathbb{R}^{n}$, the process $(Y^{s, x}, Z^{s, x})$ is the
unique solution of QBSDE
\begin{eqnarray}  \label{eqr}
Y_{t}^{s,x}&=&\psi(X_{T}^{s,x})+\int_{t}^{T}f(Y_{r}^{s,x})|Z_{r}^{s,x}|^{2}dr
\\
&& + \int_{t}^{T}G(Y_{r}^{s,x})d{\eta }_{r} -\int_{t}^{T}Z_{r}^{s,x}dB_{r}.
\nonumber
\end{eqnarray}
\end{Prop}

For an introduction to the theory of viscosity solutions, we refer the
reader to \cite{visco}. To prove the existence of viscosity solutions, we
need the following touching property, see \cite{koby}.

\begin{Lem}
Let $(\xi_{t})_{0\leq t\leq T}$ be a continuous adapted process such that
\[
d\xi_{t} = \beta_{t}dt + \alpha_{t}dW_{t},
\]
where $\beta$ and $\alpha$ are continuous adapted processes such that $\beta$%
, $|\alpha|^{2}$ are integrable. If $\xi_{t} \geq 0$ a.s. for all $t$, then
for all $t$,
\[
\mathbf{1}_{\{\xi_{t} = 0\}}\alpha_{t} = 0, ~~ \text{a.s.,}
\]
\[
\mathbf{1}_{\{\xi_{t} = 0\}}\beta_{t} \geq 0, ~~ \text{a.s.}
\]
\end{Lem}

\noindent \textbf{Proof.} First notice that
\begin{equation}
\forall \,t\in \lbrack s,T],~~u(t,X_{t}^{s,x})=Y_{t}^{s,x}.  \label{max}
\end{equation}%
This is readily seen from the Markov property of the diffusion process $X$
and from the uniqueness of the solutions of the BSDE (\ref{eqr}). Then $%
u(t,x)=Y_{t}^{t,x}$. Let $\varphi \in \mathcal{C}^{1,2}([0,T]\times \mathbb{R%
}^{n})$. Let $(t,x)$ be a local Maximum of $u-\varphi $. We suppose it
global and equal to $0$, that is
\[
\varphi (t,x)=u(t,x)~~\text{and}~~\varphi (\bar{t},\bar{x})\geq u(\bar{t},%
\bar{x})~~\text{for all}~~(\bar{t},\bar{x})\in \lbrack 0,T]\times \mathbb{R}%
^{n}.
\]%
This and equality (\ref{max}) imply that
\[
\varphi (t,X_{t}^{s,x})\geq Y_{t}^{s,x}.
\]%
We want to show that $u$ is a viscosity supersolution of (\ref{pde}).
Remember that $(Y_{\cdot }^{s,x},Z_{\cdot }^{s,x})$ satisfy%
\begin{eqnarray*}
Y_{t}^{s,x} &=&Y_{T}^{s,x}+\int_{t}^{T}f(Y_{r}^{s,x})|Z_{r}^{s,x}|^{2}dr \\
&&+\int_{t}^{T}G(Y_{r}^{s,x})\dot{\eta}_{r}dr-\int_{t}^{T}Z_{r}^{s,x}dW_{r}.
\end{eqnarray*}%
We apply Itô's formula to the process $\varphi (t,X_{t}^{s,x})$, then we
obtain
\begin{eqnarray*}
\varphi (T,X_{T}^{s,x}) &=&\varphi (t,X_{t}^{s,x})+\int_{t}^{T}\left(
\varphi _{t}+\mathcal{L}\varphi \right) (r,X_{r}^{s,x})dr \\
&&+\int_{t}^{T}\left( \varphi _{x}\sigma \right) (r,X_{r}^{s,x})dW_{r}.
\end{eqnarray*}%
As $\varphi (t,X_{t}^{s,x})\geq Y_{t}^{s,x}$, the touching property gives
for all $t$,
\[
\mathbf{1}_{\{\varphi (t,X_{t}^{s,x})=Y_{t}^{s,x}\}}\Big[\varphi
_{t}(t,X_{r}^{s,x})+\mathcal{L}\varphi (t,X_{t}^{s,x})
\]%
\begin{equation}
+f(Y_{t}^{s,x})|Z_{t}^{s,x}|^{2}+G(Y_{t}^{s,x})\dot{\eta}_{t}\Big]\geq 0,~~%
\text{a.s.,}  \label{der}
\end{equation}%
\begin{equation}
\mathbf{1}_{\{\varphi (t,X_{t}^{s,x})=Y_{t}^{s,x}\}}\Big[-Z_{t}^{s,x}+\left(
\varphi _{x}\sigma \right) (t,X_{t}^{s,x})\Big]=0,~~\text{a.s.}  \label{sig}
\end{equation}%
Or $\varphi (t,x):=\varphi (t,X_{t}^{t,x})=Y_{t}^{t,x}:=u(t,x)$ for $s=t$,
then equation (\ref{sig}) gives $Z_{t}^{t,x}=\left( \varphi _{x}\sigma
\right) (t,x)$ and equation (\ref{der}) gives the expected result.

\begin{Lem}
Let $G$ be a Lipschitz on $\mathbb{R},~\eta $ a given smooth path and $\phi $
the flow defined in (\ref{flo1}). We assume that (H4)--(H6) hold and $%
(Y^{s,x},Z^{s,x})$ is the unique solution of the QBSDE (\ref{eqr}). Then $%
u(t,x)=Y_{t}^{t,x}$ is a viscosity solution to the PDE (\ref{pde}) if and
only if $v(t,x)=\phi ^{-1}(t,u(t,x))$ is a viscosity solution of the PDE
\[
\left\{
\begin{array}{l}
\partial _{t}v(t,x)+\mathcal{L}v(t,x)+\tilde{f}(t,v(t,x),v_{x}(t,x)\sigma
(t,x))=0, \\
v(T,x)=\psi (x),~~x\in \mathbb{R},%
\end{array}%
\right.
\]%
where (in what follows the $\phi $ will be evaluated at $(t,\tilde{y})$)
\[
\tilde{f}(t,\tilde{y},\tilde{z}):=\frac{1}{\partial _{y}\phi }\left( f(\phi
)|\partial _{y}\phi \tilde{z}|^{2}+\frac{1}{2}\partial _{y}^{2}\phi |\tilde{z%
}|^{2}\right) .
\]
\end{Lem}

\noindent \textbf{Proof.} Suppose
that $u$ is a viscosity solution to the PDE (\ref{pde}). The function $u$ is
continuous, hence $v$ is also continuous. We have
\[
v(T,x)=\phi ^{-1}(T,u(T,x))=\phi ^{-1}(T,\psi (x))=\psi (x).
\]%
By Lemma \ref{lemme1} and Remark \ref{rmk1}, the BSDE$(\psi (X_{T}^{s,x}),%
\tilde{f},0,0)$ has a unique solution $\tilde{Y}_{t}^{s,x}=\phi
^{-1}(t,Y_{t}^{s,x})$ and it is connected to the following PDE
%
% By setting $\tilde{u}(t,x):=\tilde{Y}_{t}^{t,x}=\phi
%^{-1}(t,Y_{t}^{t,x})=\phi ^{-1}(t,u(t,x))=v(t,x)$ a viscosity solution to
%the following PDE
\[
\left\{
\begin{array}{l}
\partial _{t}\tilde{u}(t,x)+\mathcal{L}\tilde{u}(t,x)+\tilde{f}(t,\tilde{u}%
(t,x),\tilde{u}_{x}(t,x)\sigma (t,x))=0, \\
\tilde{u}(T,x)=\psi (x),~x\in \mathbb{R},%
\end{array}%
\right.
\]
which has a viscosity solution given by
$$\tilde{u}(t,x):=\tilde{Y}_{t}^{t,x}=\phi^{-1}(t,Y_{t}^{t,x})=\phi ^{-1}(t,u(t,x))=v(t,x).$$

For the converse, we apply Itô's formula to the process $Y_{t}^{s,x}=\phi (t,%
\tilde{Y}_{t}^{s,x})$ and use the touching property.

\begin{Lem}
\label{lemme6} Let $\gamma > p\geq 1,~ \mathbf{\eta} \in \mathcal{C}%
^{p-var}_{0}([0,T],G^{[p]}(\mathbb{R}^{d}))$ and $\Phi $ the rough flow
defined in (\ref{flo2}). Under the assumption (H3), the function
\[
\tilde{f}(t, \tilde{y}, \tilde{z}):= \frac{1}{\partial_{y}\Phi}%
\left(f(\Phi)|\partial_{y}\Phi\tilde{z}|^{2} + \frac{1}{2}%
\partial_{y}^{2}\Phi|\tilde{z}|^{2}\right).
\]
is continuous and integrable.
\end{Lem}

\noindent \textbf{Proof.} The continuity of the function $f$ is obvious. The
integrability assumption can be verified using Corollary \ref{coro}.

\subsection{Main result}

\begin{Theo}
Let $\gamma >p\geq 1,~\mathbf{\eta }\in \mathcal{C}%
_{0}^{p-var}([0,T],G^{[p]}(\mathbb{R}^{d}))$ and $\Phi $ the rough flow
defined in (\ref{flo2}). Under the assumptions (H3)--(H6), $u(t,x)=Y_{t}^{t,x}$
is a viscosity solution of the following rough PDE
\[
\left\{
\begin{array}{l}
du(t,x)+\left[ \mathcal{L}u(t,x)+f(u(t,x))|u_{x}(t,x)\sigma (t,x)|^{2}\right]
dt+G(u(t,x))d{\mathbf{\eta }}_{t}=0, \\
u(T,x)=\psi (x),~x\in \mathbb{R},%
\end{array}%
\right.
\]%
where $(Y^{s,x},Z^{s,x})$ is the unique solution of the rough BSDE$(\psi
(X_{T}^{s,x}),f(y)z^{2},G,\mathbf{\eta })$.
\end{Theo}
%\textbf{IL Y A UNE CONFUSION ENTRE }$\tilde{Y}_{t}^{t,x}\ $\textbf{\ DEFINI AVEC $\Phi $:donc on peut mettre \`a la place petit tilde un widetilde }
%
\noindent \textbf{Proof.} The function  $v(t,x):= \widetilde{Y}_{t}^{t,x}$ is
a viscosity solution to the PDE
\[
\left\{
\begin{array}{l}
\partial _{t}v(t,x)+\mathcal{L}v(t,x)+\tilde{f}(t,v(t,x),v_{x}(t,x)\sigma
(t,x))=0, \\
v(T,x)=\psi (x),~~x\in \mathbb{R},%
\end{array}%
\right.
\]%
where $\widetilde{Y}^{s,x}$ is the unique solution of the BSDE$(\psi (X^{s,x}),%
\tilde{f},0,0)$. Here $\tilde{f}$ is the function defined in Lemma \ref{lemme6} . Define
\[
Y_{t}^{s,x}:=\Phi (t,\widetilde{Y}_{t}^{s,x}).
\]%
By Corollary \ref{coro} the process $(Y^{s,x},Z^{s,x})$ is the unique
solution of the rough quadratic BSDE$(\psi (X_{T}^{s,x}),f(y)z^{2},G,\mathbf{%
\eta })$. Hence we obtain
\[
u(t,x):=Y_{t}^{t,x}=\Phi (t,v(t,x)),
\]%
and we write formally
\[
\left\{
\begin{array}{l}
du(t,x)+\left( \mathcal{L}u(t,x)+f(u(t,x))|u_{x}(t,x)\sigma
(t,x)|^{2}\right) dt+G(u(t,x))d{\mathbf{\eta }}_{t}=0, \\
u(T,x)=\psi (x),~x\in \mathbb{R}.%
\end{array}%
\right.
\]%
\noindent\textbf{Conclusion.} \\
In this paper, we have studied existence and uniqueness of  a class of quadratic BSDE with rough drivers when the terminal condition is a  square integrable random variable. We have given some examples to recover BDSDEs and also class of quadratic BSDE perturbed by a fractional Brownian motion. In the Markovian setting we have reestablished a probabilistic representation of a viscosity solution of rough PDEs by means of the nonlinear Feymann-Kac formula. \\
Remark also that the extension to appropriate multidimensional cases are straightforward.

\end{document}